\documentclass[12pt]{article}

\usepackage[margin=1.0in]{geometry}
 
\usepackage{graphicx, subcaption}      
\usepackage{amssymb, amsmath}
\usepackage{tikz}
\usetikzlibrary{patterns,angles,quotes,arrows}
\usepackage[shortlabels]{enumitem}
\usepackage{amsthm}
\newtheoremstyle{mystyle}%    % Name
  {}%                         % Space above
  {}%                         % Space below
  {}%                         % Body font
  {}%                         % Indent amount
  {\bfseries}%                % Theorem head font
  {:}%                        % Punctuation after theorem head
  { }%                        % Space after theorem head, ' ', or \newline
  {}%                         % Theorem head spec (can be left empty, meaning `normal')

\usepackage[numbers]{natbib}

\providecommand{\keywords}[1]{\textbf{\textit{Index terms---}} #1}
\usepackage{cite}

\usepackage{lipsum}

\begin{document}
\title{Reachability of Nonlinear Systems with Unknown Dynamics}
\author{Taha Shafa and
	Melkior Ornik
	\thanks{Taha Shafa and Melkior Ornik are with the Department of Aerospace Engineering and the Coordinated Science Laboratory, University of Illinois Urbana-Champaign, Urbana, USA (e-mail: tahaas2@illinois.edu, mornik@illinois.edu).}
}

\maketitle

\begin{abstract}
Determining the reachable set for a given nonlinear control system is crucial for system control and planning. However, computing such a set is impossible if the system's dynamics are not fully known. This paper is motivated by a scenario where a system suffers an adverse event mid-operation, resulting in a substantial change to the system's dynamics, rendering them largely unknown. Our objective is to conservatively approximate the system's reachable set solely from its local dynamics at a single point and the bounds on the rate of change of its dynamics. We translate this knowledge about the system dynamics into an ordinary differential inclusion. We then derive a conservative approximation of the velocities available to the system at every system state. An inclusion using this approximation can be interpreted as a control system; the trajectories of the derived control system are guaranteed to be trajectories of the unknown system. To illustrate the practical implementation and consequences of our work, we apply our algorithm to a simplified model of an unmanned aerial vehicle.

\bigskip
    
    \textbf{Notice of Previous Publication.} \normalfont{This manuscript substantially improves the work of \citep{9304326}. Theory has been generalized to include a class of non-invertible matrices and improved to provide a larger set of reachable states. All lemmas, corollaries, and Theorems 1, 2, and 4 are entirely novel. Theorem 3 has been slightly modified from existing theorems in \citep{9304326} given our new results.}
\end{abstract}


\section{Introduction}

Damage to a control system can cause significant change to its dynamics. In order to avoid endangering people located in the vicinity of the system, it is crucial to understand the system's remaining capabilities. Motivated by specific examples like an aircraft losing a wing \citep{nguyen2008flight} or a UAV becoming damaged in an urban environment \citep{chowdhary2013guidance}, \citep{jourdan2010enhancing} our goal is to conservatively approximate the unknown system's set of reachable states \citep{brockett1976nonlinear}, \citep{isidori2013nonlinear} while assuming minimal knowledge about the system dynamics. We call such a set the \textit{guaranteed reachable set (GRS)}. 

The primary contribution of this paper is to provide a meaningful underapproximation of the GRS of a control-affine system. We assume the only available information at the time of computation consists of (i) local dynamics at a single point, which can be obtained with an arbitrarily small error from applying test control inputs over a short period of time \citep{ornik2019control}, and (ii) Lipschitz bounds on the rate of change of the system's dynamics provided by prior knowledge of the system design and physical laws. The reachable set of an unknown system is impossible to compute. Prior work determining the reachable set often focused on overapproximations \citep{ornik2019control},\citep{mitchell2003overapproximating}. Without discussing reachable sets, the work in \citep{ornik2019control}, operating under similar assumptions as our paper, focused on \textit{optimistic reachability}, i.e., attempting to reach a particular objective while there exists any chance of reaching it. Conversely, this paper computes states that are guaranteed to be reachable using admissible control signals.

Apart from \citep{ornik2019control}, work on reachability under uncertainty considered computation of reachable sets with dynamics generated by a finite number of uncertain parameters \citep{filippova2018estimates}, \citep{rungger2018accurate} or having bounded disturbances \citep{mitchell2005time}, \citep{zhang2014reachable}. Work in adaptive and robust control \citep{mitchell2003overapproximating}, \citep{dullerud2013course}, \citep{ioannou2012robust} assumes more knowledge on the magnitude or structure of the system dynamics, and classical data-driven learning methods \citep{nguyen2008flight}, \citep{brunton2016discovering}, \citep{chen2018data} collect data through repeated system runs. In contrast, this paper contains substantially fewer assumptions and focuses on deriving as much information as possible for the GRS based on one system run. 

Our approach relies on the interpretation of a control system as a differential inclusion \citep{aubin2012differential}, \citep{bressan2007introduction}, \citep{smirnov2002introduction}, \citep{kurzhanski2000ellipsoidal} whose right hand side equals the set of velocities that the system can achieve at every state in the state space. While, for an unknown system, exact velocity sets may not be available anywhere, we can determine the family of all velocity sets that are consistent with prior knowledge of the local dynamics at a single point and Lipschitz bounds on the rate of change of the system's dynamics. The intersection of all elements of such a family is defined as the \textit{guaranteed velocity set (GVS)}. Such a set is difficult to express in closed form, however we can analytically derive its underapproximation to compute an underapproximation of the GRS.

The outline of the paper is as follows: we discuss the problem statement in greater detail in Section II. We then proceed to derive two simply expressible sets in Section III, with one set being a ball, and the other being a more complex convex set, later utilized to derive a polygon such that both the ball and polygon are contained within the GVS. Next, in Section IV, we use these sets to derive two classes of simple control-affine systems whose reachable sets are contained in the GRS. Lastly, we present numerical examples in Section V with a brief discussion on the implementation of our method. 

\subsection{Notation}
The set of all matrices with $n$ rows and $m$ columns is denoted by $\mathbb{R}^{n\times m}$. For any vector $v$, $\|v\|$ denotes its Euclidean norm and $\|v\|_1$ denotes its 1-norm. For any matrix $M$, $M^T$ denotes its transpose and $\|M\|$ denotes its Euclidean norm: $\|M\| =$ max\textsubscript{${\|v\|=1}$} $\|Mv\|$. Equivalently, $\|M\| = \sigma_1(M)$, where $\sigma_i(M)$ represents the $i$-th singular value of $M$. We also let $M^\dagger$ denote the Moore-Penrose pseudoinverse, $\text{Im}(M)$ denote the image (range space) of $M$, and Ker$(M)$ denote the kernel (null space) of $M$. For matrices $M$ and $N$, we will say $M \in \text{Imm}(N)$ if $M = NP$ for some matrix $P$. Notation $\mathbb{B}^n(a;b)$ denotes a closed ball in $\mathbb{R}^n$ centered at $a \in \mathbb{R}^n$ with radius $b \geq 0$ under the Euclidean norm. Set $C_L(\mathbb{R}^{n};\mathbb{R}^{m})$ denotes the set of all functions $f: \mathbb{R}^{n} \to \mathbb{R}^{m}$ with a Lipschitz constant $L$, i.e., the set of all functions $f$ that satisfy $\|f(x)-f(y)\| \leq L\|x-y\|\hspace{1 mm} \text{for all} \hspace{1 mm} x$,\hspace{1 mm}$y \in \mathbb{R}^n$. Notation $a+B\mathcal{X}$ where $a \in \mathbb{R}^n$, $B \in \mathbb{R}^{n\times m}$, and $\mathcal{X} \subseteq \mathbb{R}^m$ denotes the set $a+B\mathcal{X} = \{a + Bx~|~x\in \mathcal{X}\}$.
\section{Problem Statement}

Throughout the paper, we attempt to meaningfully underapproximate the reachable set of a nonlinear control-affine system $\mathcal{M}(f,G)$ defined by

\begin{equation}\label{eqn:PrbStatement}
    \dot{x}(t) = f(x(t)) + G(x(t))u(t), \indent x(0) = x_0,
\end{equation}

\noindent where all $t \geq 0$, $x(t) \in \mathbb{R}^n$, admissible inputs $u(t) \in \mathcal{U} = \mathbb{B}^m(0;1)$, which is a common setting in reachability analysis \citep{margaliot2007reachable}, \citep{vinter1980characterization}, and functions $f:\mathbb{R}^n \to \mathbb{R}^n$, and $G:\mathbb{R}^n \to \mathbb{R}^{n\times m}$ are globally Lipschitz continuous with Lipschitz constants $L_f \geq 0$ and $L_G \geq 0$, i.e., $f \in C_{L_f}(\mathbb{R}^n; \mathbb{R}^n)$ and $G \in C_{L_G}(\mathbb{R}^n; \mathbb{R}^{n\times m})$. Cases where $L_f = 0$ or $L_G = 0$ are simple, thus for the remainder of the paper we assume $L_f > 0$ and $L_G > 0$. Noting that any sufficiently smooth function is globally Lipschitz continuous on any compact set, the theory developed in this paper can also be applied to a system whose states are guaranteed to be bounded, which is naturally true for a large class of systems \citep{chowdhary2013guidance}, \citep{jourdan2010enhancing}. Without loss of generality, we assume $x_0 = 0$.

\subsection{Assumptions and Technical Requirements}
In order to approximate the reachable set of \eqref{eqn:PrbStatement}, the technical work of \citep{9304326} requires full actuation at $x_0 = 0$, i.e., $m=n$ with $G(0)$ being full rank. We relax these requirements, so the system $\mathcal{M}(f,G)$ need not be fully actuated at $x_0 = 0$, that is, $m$ does not necessarily equal $n$. Instead of assuming full actuation, we use the following assumption:

\textit{Assumption 1}: Functions $f$ and $G$ are of the form $f(x) = Rr(x)$ and $G(x) = RH(x)$ where $R \in \mathbb{R}^{n\times m}$, $r(x) \in \mathbb{R}^m$, and $H(x) \in \mathbb{R}^{m\times m}$ such that $H(0)$ is invertible.

The case of full actuation in \citep{9304326} corresponds the case where $R$ = $I_{n}$ in Assumption 1. Motivated by the online learning technique introduced in \citep{ornik2019control}, we make the following assumption about the knowledge regarding the system dynamics. 

\textit{Assumption 2}: Bounds $L_f$ and $L_G$ are known, as well as values $f(0)$ and $G(0)$ such that $G(0) \neq 0$. 

Note we are not assuming any knowledge about matrix $R$; we only assume such an $R$ exists. We only consider the case where $G(0) \neq 0$ because it is otherwise impossible to determine anything about guaranteed velocities at states $x\neq 0$. 

It is easily shown that $\text{Im}(G(0)) = \text{Im}(R)$ under the conditions of Assumption 1. For our future results, it is important to determine the set of $x$ such that $\text{Im}(G(0)) = \text{Im}(G(x))$. Crucially, we show that the images of $G(x)$ and $G(0)$ are equal in some neighborhood of $0$.

\textit{Lemma 1}: Under Assumptions 1 and 2, if $\|x\| < \frac{\|G(0)^{\dagger}\|^{-1}}{L_G}$, then $\text{Im}(G(x)) = \text{Im}(G(0))$.

\begin{proof}
Weyl's inequality for singular values dictates that singular values as functions on matrices are uniformly Lipschitz with respect to the operator norm \citep{stewart1998perturbation}: $\|\sigma_s(G(x)) - \sigma_s(G(0))\| \leq \|G(x) - G(0)\| \leq L_G\|x\|$, such that $1 \leq s \leq r$ with $r = \text{rank}(G(0))$. The Eckhart-Young-Mirsky theorem \citep{strang2016introduction} along with the singular value decomposition of $G(0)$ show that $\sigma_r(G(0)) = \|G(0)^\dagger\|^{-1}$ is the smallest non-zero singular value of $G(0)$. For $x$ that satisfies $\|x\| < \frac{\|G(0)^\dagger\|^{-1}}{L_G}$, we thus have $\|\sigma_s(G(x)) - \sigma_s(G(0))\| < \sigma_r(G(0))$, i.e., $\sigma_s(G(x)) > 0$. Therefore rank$(G(x)) \geq \text{rank}(G(0))$. 

We defined $G(x) = RH(x)$, so $\text{Im}(G(x)) \subset \text{Im}(R)$. Since $\text{Im}(G(0)) = \text{Im}(R)$, then $\text{Im}(G(x)) \subset \text{Im}(G(0))$. Knowing rank$(G(x)) \geq \text{rank}(G(0))$, we conclude $\text{Im}(G(x)) = \text{Im}(G(0))$.

\end{proof}

\subsection{Guaranteed Reachable Set}

Let us denote a set $\mathcal{D}_{con} \subseteq C_{L_f}(\mathbb{R}^n; \mathbb{R}^n)\times C_{L_G}(\mathbb{R}^n;\mathbb{R}^{n \times m})$ as the set of all pairs $f$ and $G$ consistent with Assumptions 1 and 2. We want to underapproximate the set of reachable states given solely the knowledge of $\mathcal{D}_{con}$. We first define the \textit{(forward) reachable set} $\mathcal{R}^{\hat{f},\hat{G}}(T, x_0) = \{\phi^{\hat{f},\hat{G}}_u(t; x_0) \hspace{1 mm} | \hspace{1 mm} u : [0,T] \to \mathcal{U}, t\in [0,T]\}$, where $\phi^{\hat{f},\hat{G}}_u(\cdot; x_0)$ denotes the controlled trajectory of the system $\mathcal{M}(\hat{f},\hat{G})$ with control signal $u$ and $\phi^{\hat{f},\hat{G}}_u(0;x_0) = 0$. 

Let $T \geq 0$. We describe the \textit{guaranteed reachable set (GRS)} as:

\begin{equation} \label{eqn:PrbStatement_5}
    \mathcal{R}^\mathcal{G}(T,x_0) =  \bigcap_{(\hat{f},\hat{G}) \hspace{1 mm} \in \hspace{1 mm} \mathcal{D}_{con}} \mathcal{R}^{\hat{f},\hat{G}}(T,x_0).
\end{equation}

\noindent The GRS describes the set of all states that are reachable by any system consistent with our knowledge of the system dynamics. 

\textit{Problem 1:} Determine or underapproximate the GRS.

To solve Problem 1, we first represent ordinary differential equations with control inputs as an \textit{ordinary differential inclusion (ODI)}. We discuss doing so in Section III. Given the assumed knowledge of the system dynamics, we develop underapproximations to the right hand side of this inclusion. In Section IV, we use these underapproximations to derive two control-affine systems whose reachable sets are subsets of $\mathcal{R}^\mathcal{G}(T,0)$.

\section{Guaranteed Velocities}

We follow the classical approach of interpreting ordinary differential equations with control inputs as inclusions \citep{aubin2012differential}, \citep{bressan2007introduction}, \citep{smirnov2002introduction}, \citep{kurzhanski2000ellipsoidal}. In this section, we formally define the \textit{guaranteed velocity set} of an unknown control system and determine analytically computable underapproximations of such a set.

\subsection{Guaranteed Velocity Set}

We define the \textit{available velocity set} of the system $\mathcal{M}(f,G)$ at state $x$ by $\mathcal{V}_x = f(x) + G(x)\mathcal{U}$, and introduce the following ODI:

\begin{equation} \label{GuaranteedVelocities1}
    \dot{x} \in \mathcal{V}_x = f(x) + G(x)\mathcal{U}, \indent x(0) = x_0.
\end{equation}

If a trajectory $\phi(\cdot ;x_0)$ satisfies \eqref{GuaranteedVelocities1}, then it obviously serves as a solution to the control system \eqref{eqn:PrbStatement} for an admissible control input, and vice versa. Given Assumption 2, set $\mathcal{V}_{x_0} = \mathcal{V}_0$ is known. The goal of this section is to provide an underapproximation for set $\mathcal{V}_x$ using sets $\mathcal{D}_{con}$ and $\mathcal{V}_0$. We first define the \textit{guaranteed velocity set (GVS)} below:

\begin{equation} \label{GuaranteedVelocities2}
    \mathcal{V}^\mathcal{G}_x = \bigcap_{(\hat{f},\hat{G}) \hspace{1 mm} \in \hspace{1 mm} \mathcal{D}_{con}}
    \hat{f}(x) + \hat{G}(x)\mathcal{U} \subseteq \mathcal{V}_x.
\end{equation}

The GVS $\mathcal{V}^\mathcal{G}_x$ is the set of all velocities that can be taken by all systems consistent with the assumed knowledge of the dynamics. Let us consider the following ODI:

\begin{equation} \label{GuaranteedVelocities3}
    \dot{x} \in \mathcal{V}^\mathcal{G}_x, \indent x_0 = 0.
\end{equation}

\noindent If $\mathcal{V}^\mathcal{G}_{\phi(T;x_0)} = \emptyset$, we will consider by convention that the trajectory of \eqref{GuaranteedVelocities3} ceases to exist at time $T$. The following proposition then holds directly from \eqref{eqn:PrbStatement_5} and \eqref{GuaranteedVelocities2}.

\textit{Proposition 1:} Let $T \geq 0$. If a trajectory $\phi$ : $[0,+\infty) \to \mathbb{R}^n$ satisfies \eqref{GuaranteedVelocities3} at all times $t \leq T$, then $\phi(T) \in \mathcal{R}^\mathcal{G}(T,x_0)$.

Proposition 1 implies that the reachable set of \eqref{GuaranteedVelocities3} is a subset of $\mathcal{R}^\mathcal{G}(T,x_0)$. As briefly described in \citep{9304326}, these sets are not necessarily equal; establishing conditions for the equality of the reachable set of \eqref{GuaranteedVelocities3} and $\mathcal{R}^\mathcal{G}(T,x_0)$ is an open problem for future work.

\subsection{Ball Underapproximation}

Proposition 1 motivates us to underapproximate $\mathcal{R}^\mathcal{G}(T,x_0)$ by determining the reachable set of \eqref{GuaranteedVelocities3}. We start by examining the geometry of the ODI. Given $x \in \mathbb{R}^n$, our previous assumptions show that 

\small
\begin{equation} \label{GuaranteedVelocities4}
\begin{gathered}
    \{(\hat{f}(x),\hat{G}(x)) \hspace{1 mm} | \hspace{1 mm} (\hat{f},\hat{G}) \in \mathcal{D}_{con}\} = \\
    (\mathbb{B}^n(f(0);L_f\|x\|) \cap \text{Im}(R)) \times (\mathbb{B}^{n \times m}(G(0);L_G\|x\|) \cap \text{Imm}(R)).
\end{gathered}
\end{equation}
\normalsize

Given that $\mathcal{U} = \mathbb{B}^m(0;1)$, $\mathcal{V}^\mathcal{G}_x$ is an intersection of infinitely many ellipsoids $a + B\mathcal{U}$, where $a \in \mathbb{B}^n(f(0);L_f\|x\|) \cap \text{Im}(R)$ and $B \in \mathbb{B}^{n \times m}(G(0);L_G\|x\|) \cap \text{Imm}(R)$. An intersection of infinitely many ellipsoids is generally not a geometrically simple object \citep{boyd2004convex}. Thus, we will determine an underapproximation of $\mathcal{V}^\mathcal{G}_x$. Our approach will be to implicitly exploit convexity of $\mathcal{V}^\mathcal{G}_x$ and underapproximate the distance of the boundary of set $\mathcal{V}^\mathcal{G}_x$ from $0$ in every direction. In Theorem 1, we calculate one such underapproximation.

\textit{Theorem 1}: Let $\mathcal{U}$, $L_f$, and $L_G$ be defined as above. Let $x \in \mathbb{R}^n$ satisfy $(L_f + L_G)\|x\| < \|G(0)^{\dagger}\|^{-1}$. Define

\begin{equation} \label{eqn:Theorem_1}
    \fontsize{9.5}{5}
    \bar{\mathcal{V}}^\mathcal{G}_{x} = \mathbb{B}^{n}(f(0); \|G(0)^{\dagger}\|^{-1} - L_f\|x\| - L_G\|x\|) \hspace{1 mm} \cap \hspace{1 mm} \text{Im}(G(0)).
\end{equation}

\noindent Then, $\bar{\mathcal{V}}^\mathcal{G}_{x} \subseteq \mathcal{V}^\mathcal{G}_{x}$.

\begin{proof}
Set $\cap_{(\hat{f},\hat{G}) \in \mathcal{D}_{con}}\hat{f}(x) - f(0) + \hat{G}(x)\mathcal{U} \text{~equals~} \mathcal{V}^\mathcal{G}_x - f(0)$. On the other hand, $\cap_{(\hat{f},\hat{G}) \in \mathcal{D}_{con}} \hat{f}(x) - f(0) + \hat{G}(x)\mathcal{U} = \cap_{(\Tilde{f},\hat{G})\in \tilde{\mathcal{D}}_{con}}\Tilde{f}(x) + \hat{G}(x)\mathcal{U}$, where $\tilde{\mathcal{D}}_{con}$ is defined the same as before, just with the assumption that $f(0) = 0$. Thus, we can assume without loss of generality that $f(0) = 0$.

Let $d \in \text{Im}(G(0)) = \text{Im}(R)$ such that $\|d\| = 1$. We will prove that if $|k| \leq \|G(0)^{\dagger}\|^{-1} - L_f\|x\| - L_G\|x\|$, then equation

\begin{equation} \label{eqn:Thm1_ProblemFormulation}
    k \cdot d = \hat{f}(x) + \hat{G}(x)u,
\end{equation}

\noindent where $(\hat{f},\hat{G}) \in \mathcal{D}_{con}$, admits a solution $u \in \mathcal{U}=\mathbb{B}^{m}(0;1)$.

We subtract $\hat{f}(x)$ from both sides of \eqref{eqn:Thm1_ProblemFormulation}. Since $\hat{f}(x) \in \text{Im}(R)$ by Assumption 1 and $kd \in \text{Im}(R)$ by definition, then $kd - \hat{f}(x) \in \text{Im}(R)$. Also, $\text{Im}(\hat{G}(x)) = \text{Im}(R)$ by Lemma~1. Hence, there exists a vector $\bar{u} \in \mathbb{R}^m$ such that $kd - \hat{f}(x) = \hat{G}(x)\bar{u}$. Now, through the rank-nullity theorem \citep{strang2016introduction}, we can write $\bar{u} = u + u_2$ where $u \in \text{Im}(\hat{G}(x)^T)$ and $u_2 \in \text{Ker}(\hat{G}(x))$. Thus, $\hat{G}(x)\bar{u} = \hat{G}(x)(u + u_2) = \hat{G}(x)u$; hence, $kd - \hat{f}(x) = \hat{G}(x)u$.

We multiply both sides of $kd - \hat{f}(x) = \hat{G}(x)u$ on the left by $\hat{G}(x)^\dagger$, resulting in $\hat{G}(x)^\dagger (kd - \hat{f}(x)) = \hat{G}(x)^\dagger \hat{G}(x)u$. The term $\hat{G}(x)^\dagger \hat{G}(x)u$ results in the projection of $u$ onto the $\text{Im}(\hat{G}(x)^T)$ \citep{strang2016introduction}. Given that $u \in \text{Im}(\hat{G}(x)^T)$, by definition of a projection, $\hat{G}(x)^\dagger (kd - \hat{f}(x)) = \hat{G}(x)^\dagger \hat{G}(x)u = u$. Thus, if we prove that:

\begin{equation}
    \|\hat{G}(x)^{\dagger}(k \cdot d-\hat{f}(x))\| \leq 1,
    \label{eqn:thm1_1}
\end{equation}

\noindent we will have $\|u\| \leq 1$, i.e., $u \in \mathcal{U}$. Utilizing $\|d\| = 1$ along with the product and triangle inequalities for matrices, we arrive at \eqref{eqn:thm1_4} and \eqref{eqn:thm1_5}:

\begin{align}
    \|\hat{G}(x)^{\dagger}(k \cdot d-\hat{f}(x))\| \leq |k|\|\hat{G}(x)^{\dagger}~ d\| + \|\hat{G}(x)^{\dagger} \hat{f}(x)\| \label{eqn:thm1_4}, \\
    \leq |k|\|\hat{G}(x)^{\dagger}\| + \|\hat{G}(x)^{\dagger}\| \|\hat{f}(x)\|. \label{eqn:thm1_5}
\end{align}

From \eqref{eqn:thm1_5} it follows that the set of all $k$ that satisfy $|k| \|\hat{G}(x)^{\dagger}\| + \|\hat{G}(x)^{\dagger}\| \|\hat{f}(x)\| \leq 1$ is a subset of all $k$ that satisfy \eqref{eqn:thm1_1}. In other words, if:

\begin{equation}
    |k| \leq \|\hat{G}(x)^\dagger\|^{-1} - \|\hat{f}(x)\|,
    \label{eqn:thm1_2}
\end{equation}

\noindent then $k$ satisfies \eqref{eqn:thm1_1}. We note $\|\hat{G}(x)^\dagger\| \neq 0$ from the definition of the Moore-Penrose pseudoinverse and because $G(x) \neq 0$ from Lemma~1. 

By Weyl's inequality for singular values \citep{stewart1998perturbation} and Assumption 2, we obtain the following inequalities: $\|\hat{G}(x)^{\dagger}\|^{-1} \geq \|G(0)^{\dagger}\|^{-1} - L_G\|x\|$ and $\|\hat{f}(x)\| \leq L_f\|x\|$. Thus, since we assumed that $k$ satisfies:

\begin{equation}
    |k| \leq \|G(0)^{\dagger}\|^{-1} - L_f\|x\| - L_G\|x\|,
    \label{eqn:thm1_6}
\end{equation}

\noindent it satisfies \eqref{eqn:Thm1_ProblemFormulation}.

\end{proof}

We slightly generalize Theorem 1 by also considering all $x$ which satisfy $(L_f + L_G)\|x\| = \|G(0)^\dagger\|^{-1}$.

\textit{Corollary 1}: If $\|x\| = \frac{\|G(0)^\dagger\|^{-1}}{L_f + L_G}$, then $f(0) \in \mathcal{V}^\mathcal{G}_x$.

\begin{proof}
Let us take $x$ such that $\|x\| = \frac{\|G(0)^\dagger\|^{-1}}{L_f + L_G}$ and a sequence $x_1, x_2,...$ such that $\|x_i\| < \frac{\|G(0)^\dagger\|^{-1}}{L_f + L_G}$ and $\|x_i\| \to x$ as $i \to \infty$. Theorem 1 shows that for all $i$, there exists a $u_i \in \mathcal{U}$ such that $f(0) = \hat{f}(x_i) + \hat{G}(x_i)u_i$. Thus, because $\mathcal{U}$ is a compact set, there exists a subsequence $u_{p_1}, u_{p_2},...$ which converges to some $u^* \in \mathcal{U}$ \citep{ross2013elementary}. Since $x_{p_i} \to x$ and $u_{p_i} \to u^*$, we have $f(0) = \hat{f}(x_{p_i}) + \hat{G}(x_{p_i})u_{p_i} \to \hat{f}(x) + \hat{G}(x)u^*$. Therefore, $f(0) = \hat{f}(x) + \hat{G}(x)u^*$.
    
\end{proof}

\subsection{Advanced Convex Underapproximation}

For all $x$ that satisfy $\|x\| \leq \frac{\|G(0)^\dagger\|^{-1}}{L_f + L_G}$, we now have a set $\bar{\mathcal{V}}^\mathcal{G}_x$ guaranteed to be a subset of $\mathcal{V}^\mathcal{G}_x$. Such a set is a projection of a ball onto $\text{Im}(G(0))$. However, there are instances, particularly when singular values of $G(0)$ are far apart, where a ball may be a poor underapproximation of $\mathcal{V}^\mathcal{G}_x$. Consequently, we derive a new underapproximated set $\bar{\bar{\mathcal{V}}}^\mathcal{G}_x$ in Theorem 2.

\textit{Theorem 2}: Let $\mathcal{U}$, $L_f$, and $L_G$ be defined as above. Let $\mu = 1$ if rank($G(0)) = m = n$, $\mu = \sqrt{2}$ if rank($G(0)) = \text{min}(m,n)$ and $m \neq n$, $\mu = \frac{1 + \sqrt{5}}{2}$ if rank($G(0)) < \text{min}(m,n)$, and let $x$ satisfy $(L_f + L_G)\|x\| \leq \|G(0)^{\dagger}\|^{-1}$. If

\small
\begin{equation}\label{eqn:Theorem_2}
\begin{gathered}
    \bar{\bar{\mathcal{V}}}^\mathcal{G}_x = \{f(0) + kd \hspace{1 mm} | \hspace{1 mm} \|d\| = 1, d\in \text{Im}(R), \hspace{1 mm} 0 \leq k \leq K(d)\} \\
    \text{s.t.} \hspace{1 mm} K(d) = \frac{\|G(0)^\dagger\|^{-1} - L_G\|x\| - L_f\|x\|}{\|G(0)^\dagger d\|(\|G(0)^\dagger\|^{-1} - L_G\|x\|) + \mu\|G(0)^\dagger\| L_G\|x\|},  
\end{gathered}
\end{equation}

\normalsize
\noindent then $\bar{\bar{\mathcal{V}}}^\mathcal{G}_x \subseteq \mathcal{V}^\mathcal{G}_x$.

\begin{proof}
As in Theorem 1, we will show that for $k$ and $d$ given in \eqref{eqn:Theorem_2}, equation \eqref{eqn:Thm1_ProblemFormulation} admits a solution $u \in \mathcal{U}$. Like in the proof of Theorem 1, without loss of generality, we set $f(0) = 0$. From inequality \eqref{eqn:thm1_4}, it follows that the set of all $k$ that satisfy 

\begin{equation}
    |k| \|\hat{G}(x)^{\dagger}d\| + \|\hat{G}(x)^{\dagger}\| \|\hat{f}(x)\| \leq 1
    \label{eqn:thm2_1}
\end{equation}

\noindent is a subset of all $k$ that satisfy \eqref{eqn:thm1_1}. The term $\|G(x)^\dagger\|$ is the inverse of the smallest nonzero singular value of $G(x)$. Therefore, using Weyl's inequality for singular values \citep{stewart1998perturbation}, we have $\|\hat{G}(x)^\dagger\| \leq (\|G(0)^\dagger\|^{-1} - L_G\|x\|)^{-1}$. By Assumption~2, we also have $\|\hat{f}(x)\| \leq L_f\|x\|$. We conclude that the set of all $k$ that satisfy

\begin{equation} \label{eqn:thm2_2}
    |k|\|\hat{G}(x)^{\dagger}d\| + \frac{L_f\|x\|}{\|G(0)^{\dagger}\|^{-1} - L_G\|x\|} \leq 1
\end{equation}

\noindent is a subset of all $k$ that satisfy \eqref{eqn:thm1_1}.

Next, we bound $\|\hat{G}(x)^{\dagger}d\|$ using the product and triangle inequalities for matrices:

\begin{equation}\label{eqn:thm2_3}
\begin{gathered}
\begin{aligned}
    \|\hat{G}(x)^{\dagger}d\| \leq \|G(0)^{\dagger}d\| + \|(\hat{G}(x)^\dagger - G(0)^\dagger)d\| \\
    \leq \|G(0)^{\dagger}d\| + \|\hat{G}(x)^\dagger - G(0)^\dagger\|.
\end{aligned}
\end{gathered}
\end{equation}

\noindent By Lemma 1, $\text{Im}(G(0)) = \text{Im}(\hat{G}(x))$. Therefore, we can apply the inequality $\|\hat{G}(x)^\dagger - G(0)^\dagger\| \leq \mu \|\hat{G}(x)^\dagger\| \|G(0)^\dagger\| \|\hat{G}(x) - G(0)\|$ (Theorem 3.3 in \citep{stewart1977perturbation}). We can now rewrite the upper bound on $\|\hat{G}(x)^{\dagger}d\|$:

\begin{equation} \label{eqn:thm2_4}
    \|\hat{G}(x)^{\dagger}d\| \leq \|G(0)^{\dagger}d\| + \mu\|\hat{G}(x)^\dagger\| \|G(0)^\dagger\| \|\hat{G}(x) - G(0)\|.
\end{equation}

We again use $\|\hat{G}(x)^\dagger\| \leq (\|G(0)^\dagger\|^{-1} - L_G\|x\|)^{-1}$. According to Weyl's inequality for singular values, $\|\hat{G}(x) - G(0)\| \leq L_G\|x\|$. Therefore, all $\|\hat{G}(x)^{\dagger}d\|$ that satisfy \eqref{eqn:thm2_4} will also satisfy

\begin{equation} \label{eqn:thm2_5}
    \|\hat{G}(x)^{\dagger}d\| \leq \|G(0)^{\dagger}d\| + \frac{\mu\|G(0)^\dagger\| L_G\|x\|}{\|G(0)^\dagger\|^{-1} - L_G\|x\|}.
\end{equation}

By plugging \eqref{eqn:thm2_5} into \eqref{eqn:thm2_2}, we obtain all $k$ that satisfy 

\small
\begin{equation*}
\begin{split}
    |k| \hspace{1 mm} \left(\|G(0)^{\dagger}d\| + \frac{\mu\|G(0)^\dagger\| L_G\|x\|}{\|G(0)^\dagger\|^{-1} - L_G\|x\|}\right) \\ +
     \frac{L_f\|x\|}{\|G(0)^\dagger\|^{-1} - L_G\|x\|} \leq 1,
\end{split}
\end{equation*}
\normalsize

\noindent i.e.,

\small
\begin{equation} \label{eqn:thm2_6}
\begin{gathered}
    |k| \leq \frac{\|G(0)^\dagger\|^{-1} - L_G\|x\| - L_f\|x\|}{\|G(0)^{\dagger}d\|(\|G(0)^\dagger\|^{-1} - L_G\|x\|) + \mu\|G(0)^\dagger\| L_G\|x\|},  
\end{gathered}
\end{equation}
\normalsize

\noindent also satisfy \eqref{eqn:thm1_1}.
\end{proof}

By rewriting the denominator of $K(d)$ in Theorem 2 as $\|G(0)^\dagger\|^{-1}\|G(0)^\dagger d\| + (\mu\|G(0)^\dagger\| - \|G(0)^\dagger d\|)L_G\|x\|$, we can easily see that it cannot be negative since $\mu \geq 1$ and $\|G(0)^\dagger\| \geq \|G(0)^\dagger d\|$. We now compare the derived sets $\bar{\mathcal{V}}^\mathcal{G}_x$ and $\bar{\bar{\mathcal{V}}}^\mathcal{G}_x$ in Corollary 2.

\textit{Corollary 2}: For invertible matrices $G(0)$, $\bar{\mathcal{V}}^\mathcal{G}_x \subseteq \bar{\bar{\mathcal{V}}}^\mathcal{G}_x$.

\begin{proof}
Given the numerator of $K(d)$ in \eqref{eqn:Theorem_2} is identical to the radius of $\bar{\mathcal{V}}^\mathcal{G}_x$ derived in Theorem 1, $\bar{\mathcal{V}}^\mathcal{G}_x$ is contained within $\bar{\bar{\mathcal{V}}}^\mathcal{G}_x$ if and only if

\begin{equation} \label{Thm2_denominator}
    \|G(0)^{\dagger}d\|(\|G(0)^\dagger\|^{-1} - L_G\|x\|) + \mu\|G(0)^\dagger\| L_G\|x\| \leq 1.
\end{equation}

\noindent Inequality \eqref{Thm2_denominator} is obviously equivalent to 

\begin{equation} \label{eqn:thm2_7}
    \|G(0)^\dagger d\|(\|G(0)^\dagger\|^{-1} - L_G\|x\|) \leq 1 - \mu\|G(0)^\dagger\|L_G\|x\|.
\end{equation}

\noindent Dividing by $\|G(0)^\dagger\|^{-1} - L_G\|x\| > 0$ results in

\begin{equation} \label{eqn:thm2_8}
\begin{gathered}
    \|G(0)^\dagger d\| \leq \frac{\|G(0)^\dagger\|(1 - \mu\|G(0)^\dagger\|L_G\|x\|)}{1 - \|G(0)^\dagger\|L_G\|x\|}.
\end{gathered}
\end{equation}

We see in the case of invertible $G(0)$ where $\mu = 1$, \eqref{eqn:thm2_8} reduces to $\|G(0)^\dagger d\| \leq \|G(0)^\dagger\|$, which holds true for all d such that $\|d\| = 1$. Therefore, for invertible $G(0)$, $\bar{\mathcal{V}}^\mathcal{G}_x$ is contained within $\bar{\bar{\mathcal{V}}}^\mathcal{G}_x$. 

\end{proof}

\begin{figure}[htbp!]
	\centering
	\includegraphics[scale=0.5]{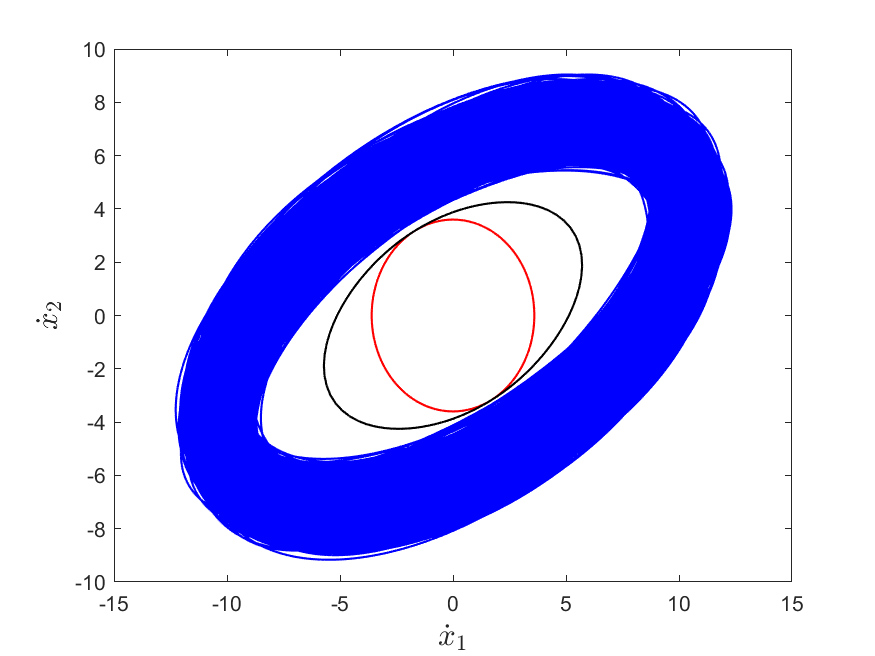}
	\caption{Velocity sets plotted for a system with $f(0) = \begin{bmatrix}  0 \\ 0  \end{bmatrix}$, $G(0) = \begin{bmatrix} 10 & 3 \\ 2 & 7 \end{bmatrix}$, $\|x\| = 1$, and $L_f = L_G = 1$. The blue curves represent the boundaries of the available velocity sets for all systems such that $(\hat{f}, \hat{G}) \in \mathcal{D}_{con}$. The intersection of these sets produces the guaranteed velocity set $\mathcal{V}^\mathcal{G}_x$ (white). Underapproximations $\bar{\mathcal{V}}^\mathcal{G}_x$ (bounded in red) and $\bar{\bar{\mathcal{V}}}^\mathcal{G}_x$ (bounded in black) are both contained in $\mathcal{V}^\mathcal{G}_x$.}
    \label{fig:AcademicSystem_VelocitySets}
\end{figure}

Our calculations result in two derived sets. Corollary 2 proves $\bar{\mathcal{V}}^\mathcal{G}_x \subseteq \bar{\bar{\mathcal{V}}}^\mathcal{G}_x$ when $G(0)$ is invertible. Hence, Theorem 2 provides a better underapproximation than the one obtained in \citep{9304326}, under the assumptions present in that paper. Figure \ref{fig:AcademicSystem_VelocitySets} illustrates a simple example of the GVS $\mathcal{V}^\mathcal{G}_x$ and its approximations obtained by Theorems 1 and 2. Consistent with the results of Corollary 2, since $G(0)$ is invertible, $\bar{\mathcal{V}}^\mathcal{G}_x \subseteq \bar{\bar{\mathcal{V}}}^\mathcal{G}_x$. 

For general cases, $\bar{\mathcal{V}}^\mathcal{G}_x \not\subseteq \bar{\bar{\mathcal{V}}}^\mathcal{G}_x$. Thus, we can take the union of both of sets to determine a larger set of guaranteed velocities. We denote such a set by

\begin{equation*}
    \hat{\mathcal{V}}^\mathcal{G}_x = \bar{\mathcal{V}}^\mathcal{G}_x \cup \bar{\bar{\mathcal{V}}}^\mathcal{G}_x.
\end{equation*}

\noindent In Section IV, we show how  $\hat{\mathcal{V}}^\mathcal{G}_x$ can be used to identify a polygon that can generate a control system with solutions that satisfy \eqref{eqn:PrbStatement}.

\section{Reachable Set}

In this section, we aim to utilize the derived sets $\bar{\mathcal{V}}^\mathcal{G}_x$ and $\hat{\mathcal{V}}^\mathcal{G}_x$ to determine a set of trajectories guaranteed to satisfy \eqref{eqn:PrbStatement}. First, we use $\hat{\mathcal{V}}^\mathcal{G}_x$ to identify a polygon contained within $\mathcal{V}^\mathcal{G}_x$. 

\textit{Lemma 2}: Let $\mathcal{S}(x)$ be any finite set of points on the boundary of $\hat{\mathcal{V}}^\mathcal{G}_x$. Let $P(\mathcal{S}(x))$ be a convex hull of $\mathcal{S}(x)$. Every solution to $\dot{x} \in P(\mathcal{S}(x))$ is a solution to \eqref{eqn:PrbStatement}.

\begin{proof}
    The guaranteed velocity set is the intersection of an infinitely many ellipsoids; such a set is convex \citep{boyd2004convex}. Given $\hat{\mathcal{V}}^\mathcal{G}_x \subseteq \mathcal{V}^\mathcal{G}_x$, we know $P(\mathcal{S}(x)) \subseteq \mathcal{V}^\mathcal{G}_x$.
\end{proof}

While Lemma 2 permits us to chose any set of points on the boundary of $\hat{\mathcal{V}}^\mathcal{G}_x$ as $\mathcal{S}(x)$, in the remainder of the paper, we will choose the points along the singular vectors of $G(0)$. We do so by noting from inequality \eqref{eqn:thm2_6} that the smaller the magnitude of $\|G(0)^\dagger d\|$, the larger the corresponding $K(d)$. The smallest magnitude of $\|G(0)^\dagger d\|$ for $\|d\| = 1$ will be obtained when $d$ is the singular vector corresponding to the largest singular value of $G(0)$ \citep{strang2016introduction}. We choose other points to be along other singular vectors because of the orthogonality of singular vectors.

Theorem 1 shows that the reachable set of

\begin{equation} \label{reachSet1}
    \dot{x} \in \bar{\mathcal{V}}^\mathcal{G}_x, \indent x(0) = x_0,
\end{equation}

\noindent is a subset of $\mathcal{R}^\mathcal{G}_x(T,x_0)$, while Lemma 2 shows the same for the reachable set of

\begin{equation} \label{reachSet2}
    \dot{x} \in P(\mathcal{S}(x)), \indent x(0) = x_0.
\end{equation}

\subsection{Underapproximated Control System -- Ball}

Analogous to the interpretation of dynamics \eqref{eqn:PrbStatement} as an ODI \eqref{GuaranteedVelocities1}, inclusion \eqref{reachSet1} can be interpreted as a control system

\begin{equation}\label{ODE1}
    \dot{x} = a + g(\|x\|)u, \indent x(0) = x_0,
\end{equation}

\noindent on $\{x~|~\|x\| \leq \|G(0)^\dagger\|/(L_f + L_G)\}$, with $a = f(0)$, $u \in \mathcal{U} = \mathbb{B}^m(0;1) \cap \text{Im}(G(0)^T)$, and where $g(s) = \|G(0)^\dagger\|^{-1} - (L_f + L_G)s$ if $s \leq \|G(0)^\dagger\|^{-1}/(L_G + L_f)$. We thus obtain the following result.

\textit{Theorem 3}: Let $\bar{\mathcal{R}}(T,x_0)$ be defined as the reachable set of \eqref{ODE1} at time $T$. Then, $\bar{\mathcal{R}}(T,x_0) \subseteq \mathcal{R}^\mathcal{G}(T,x_0)$.

\begin{proof}
Proposition 1 and Theorem 1 show that $\bar{\mathcal{R}}(T,x_0) \subseteq \mathcal{R}^\mathcal{G}(T,x_0)$.
\end{proof}

We can expand on Theorem 3 using control system \eqref{ODE1} to help determine the geometric structure of $\bar{\mathcal{R}}(T,x_0)$.

\textit{Corollary 3}: Consider a control system of the form $\dot{x} = (b - c\|x\|)u$ defined on some ball $\mathbb{B} \subseteq \mathbb{R}^n$, such that $x \in \mathbb{R}^n$, $b,c \in \mathbb{R}$, and $u \in \mathcal{U}$. Then, the reachable set $\bar{\mathcal{R}}(T,x_0)$ is a ball in $\mathbb{B}$.

\begin{proof}
    Let $z = Rx$ such that $R \in \mathbb{R}^{n\times n}$ is any orthonormal matrix. Obviously, $\dot{z} = R\dot{x} = (b - c\|x\|)Ru$. Note that $\|Rx\| = \|x\| = \|z\|$, thus we have $\dot{z} = (b + c\|z\|)Ru$. Similarly, let $Ru = v$. Then, $\|v\| = \|Ru\| = \|u\|$. Therefore, $\dot{z} = (b - c\|z\|)v$ such that $v \in \mathcal{U}$, so the reachable set of $\dot{x} = (b - c\|x\|)u$ is invariant to all rotations. Any rotation of any point on any trajectory of the original system is thus itself on some other trajectory of the original system. Hence, as the trajectories are continuous, the reachable set of the original system is a ball.
\end{proof}

\subsection{Underapproximated Control System -- Polygon}

Analogous to Theorem 3, we can define $\hat{\mathcal{R}}(T,x_0)$ as the reachable set of \eqref{reachSet2} at time $T$, and again $\hat{\mathcal{R}}(T,x_0) \subseteq \mathcal{R}^\mathcal{G}(T,x_0)$. We follow the method above and interpret inclusion \eqref{reachSet2} as a control system defined in Theorem 4. 

\textit{Theorem 4}: Let $s \in \mathbb{R}$ and $g(s) = \|G(0)^\dagger\|^{-1} - (L_G + L_f)s$, $\alpha(s) = \|G(0)^\dagger\|^{-1} - L_Gs$, $\beta(s) = \mu\|G(0)^\dagger\|L_Gs$ with $\mu$ as defined in Theorem 2. Let $U\Sigma V^T$ be the singular value decomposition of $G(0)$ where $U = [\eta_1,...,\eta_n]$. Let $r = \text{rank}(G(0))$; we define $\Lambda(s) = \text{diag}(\lambda_i(s))$ such that $\lambda_i(s) = \text{max}\{\frac{g(s)}{\alpha(s)\|G(0)^\dagger \eta_i\| + \beta(s)},g(s),0\}$ for $i = 1,...,r$ and $\lambda_i(s) = 0$ elsewhere. 

The reachable set of \eqref{reachSet2} equals the reachable set of the control system

\begin{equation}\label{ODE2}
    \dot{x} = a + U\Lambda(\|x\|)u, \indent x(0) = x_0,
\end{equation}

\noindent on $\{x~|~\|x\| \leq \|G(0)^\dagger\|/(L_f + L_G)\}$, with $a = f(0)$ and $u \in \{u~|~\|u\|_1 \leq 1\} \cap \text{Im}(G(0)^T)$. If $\hat{\mathcal{R}}(T,x_0)$ denotes the reachable set of \eqref{ODE2}, then $\hat{\mathcal{R}}(T,x_0) \subseteq \mathcal{R}^\mathcal{G}(T,x_0)$.

\begin{proof}
    Let $\|x\|\leq \|G(0)^\dagger\|/(L_f + L_G)\}$. By Theorem 1, the boundary of $\bar{\mathcal{V}}^\mathcal{G}_x$ lies at a magnitude of $g(\|x\|)$ along any direction in $\text{Im}(G(0))$ from $f(0)$. By Theorem 2, the boundary of $\bar{\bar{\mathcal{V}}}^\mathcal{G}_x$ lies at a magnitude of $g(\|x\|)/(\alpha(\|x\|)\|G(0)^\dagger \eta_i\| + \beta(\|x\|))$ along the direction $\eta_i$ from $f(0)$, for all $i\leq r$. Since $\hat{\mathcal{V}}^\mathcal{G}_x=\bar{\mathcal{V}}^\mathcal{G}_x\cup\bar{\bar{\mathcal{V}}}^\mathcal{G}_x$, the vertex of $\mathcal{S}(x)$ in the direction $\eta_i$ thus has a magnitude of $\lambda_i(\|x\|)$, for $i\leq r$. In the direction of $\eta_i$ for $i>r$, the distance of boundary points for both $\bar{\mathcal{V}}^\mathcal{G}_x$ and $\bar{\bar{\mathcal{V}}}^\mathcal{G}_x$ from $f(0)$ equals $0$, as both of these sets are contained in $\text{Im(G(0))}$. By its construction in Lemma 2, polygon $P(\mathcal{S}(x))$ is thus given by $f(0)+U\Lambda(\|x\|)Q$, where $Q=\{u~|~\|u\|_1\leq 1\}$.
\end{proof}

We note that the dynamics of \eqref{ODE1} and \eqref{ODE2} are \textit{entirely known}. Finding $\bar{\mathcal{R}}(T,x_0)$ and $\hat{\mathcal{R}}(T,x_0)$ becomes a standard problem of determining the reachable set of a nonlinear control system \citep{sontag1988controllability}. In order to exploit previous work on computing reachable sets in $\mathbb{R}^n$, we can follow methods outlined in \citep{9304326}, where \eqref{reachSet1} and \eqref{reachSet2} can be continuously extended to all $\mathbb{R}^n$ by defining $\bar{\mathcal{V}}^\mathcal{G}_x = \{f(0)\}$ and $P(\mathcal{S}(x)) = \{f(0)\}$ for all $x$ such that $\|x\| > \frac{\|G(0)^\dagger\|^{-1}}{L_f + L_G}$. Although there is no method to analytically determine the exact reachable set for all nonlinear systems, existing level set methods find the reachable set by determining the viscosity solution to Hamilton-Jacobi equations \citep{mitchell2005time}, \citep{chen2017high}. Additional methods create an overapproximation of the true reachable set by utilizing trajectory piecewise linearized models \citep{han2006reachability} or set propagation techniques \citep{althoff2008reachability}, \citep{chachuat2015set}, \citep{ramdani2011computing}. For simplicity, in our numerical examples, we approximate the true reachable set using a Monte Carlo method by solving ODEs \eqref{ODE1} and \eqref{ODE2} with random time-varying inputs.

\section{Numerical Examples}

We consider two examples. The first shows novel theory applied to an academic three-dimensional nonlinear system with initial conditions similar to the system briefly discussed in Section III, shown in Figure \ref{fig:AcademicSystem_VelocitySets}. This illustrates an example where Theorem 1 results in a poorer underapproximaiton of the true reachable set, while utilizing $\hat{\mathcal{R}}(T,x_0)$ generates a significantly better underapproximation. The second example is a control system with decoupled quadrocoptor dynamics motivated by the scenario of landing a damaged UAV safely \citep{chowdhary2013guidance}, \citep{jourdan2010enhancing}. The goal is to determine a reachable set of pitch and roll velocities in order to help stabilize the UAV for landing. 

As mentioned at the end of Section IV, we calculate the reachable sets by numerically solving all ordinary differential equations using the standard ode45 function in MATLAB which implements a Runge-Kutta method \citep{dormand1980family} with a variable time step for efficient computation. Doing so avoids limitations of current numerical solvers such as CORA \citep{althoff2018implementation} that often rely on set propagation methods to calculate the reachable set and face issues as $\|x\|$ is not differentiable at $0$. Numerical results are supported through analytical means and theoretical results derived in Sections II, III, and IV.

\subsection{Three-Dimensional Nonlinear System}

\begin{figure*}[h]
\centering
\subfloat{%
       \includegraphics[width=0.32\linewidth]{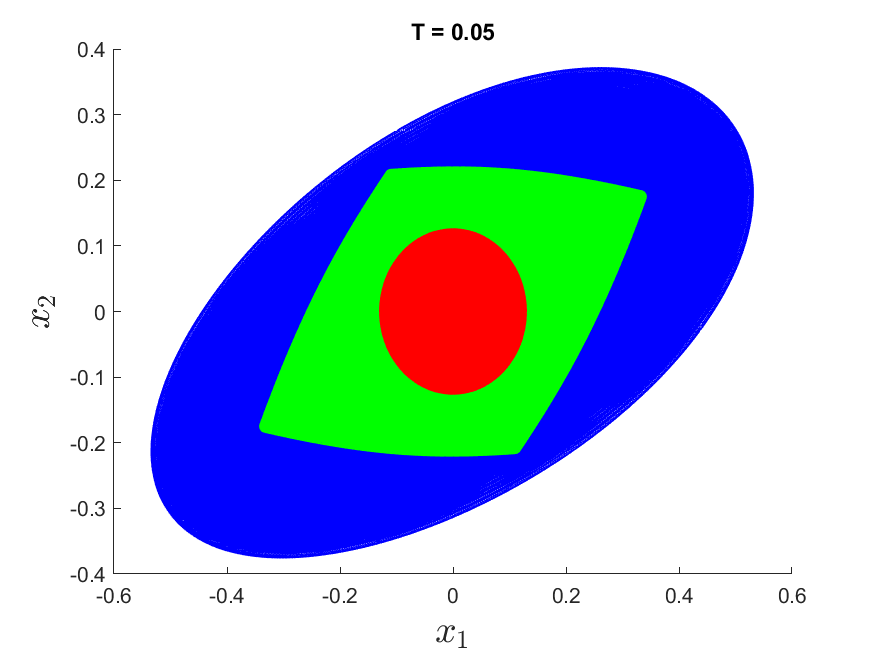}}
    \hfill
\captionsetup{justification=centering,margin=2cm}
\subfloat{%
       \includegraphics[width=0.32\linewidth]{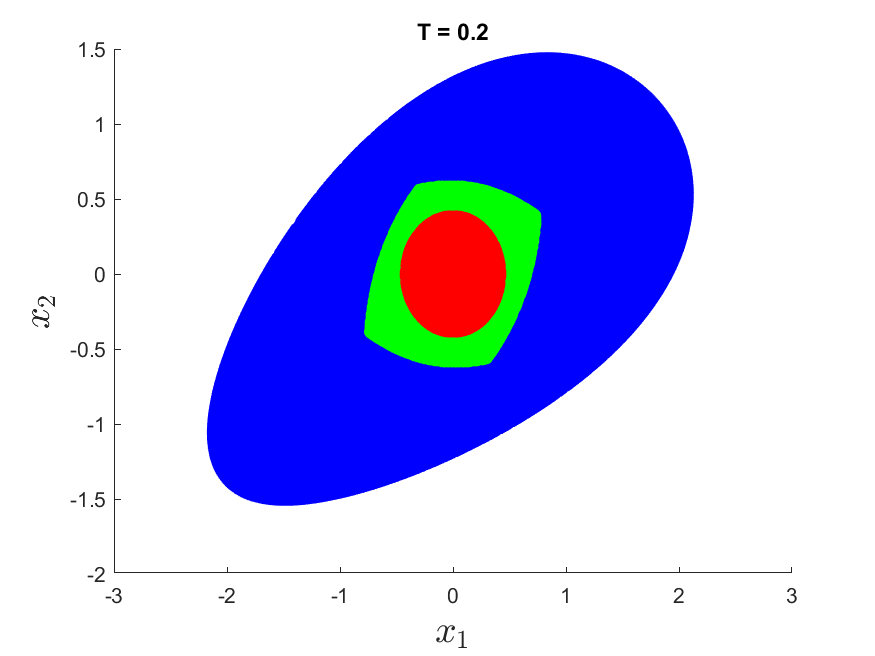}}
    \hfill
\captionsetup{justification=centering,margin=2cm}
\subfloat{%
       \includegraphics[width=0.32\linewidth]{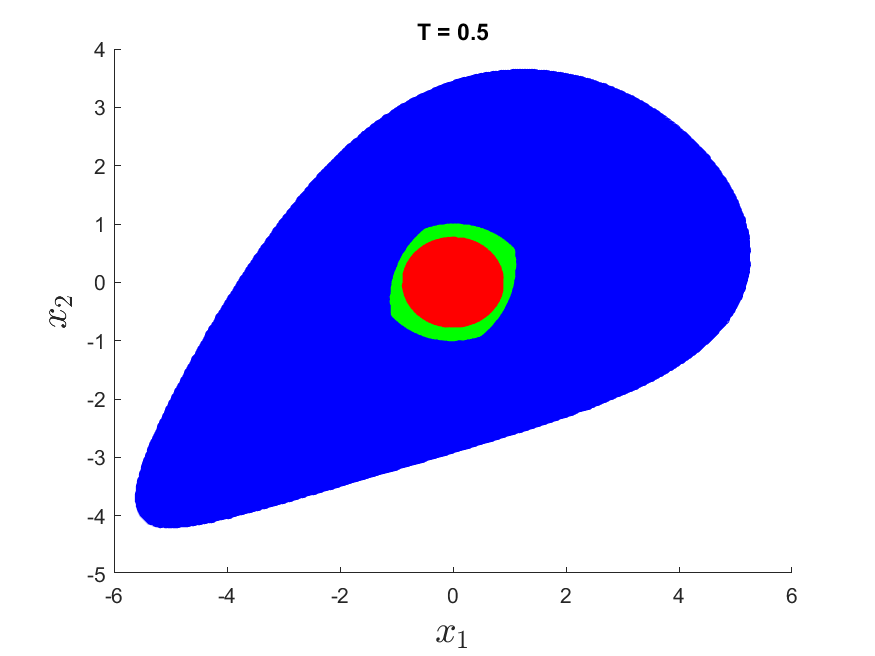}}
    \hfill
\captionsetup{justification=centering,margin=2cm}
\caption[font=small]{True reachable set (blue) with the underapproximations $\bar{\mathcal{R}}(T,0)$ (red) and $\hat{\mathcal{R}}(T,0)$ (green) numerically calculated for $T \in \{0.05, 0.2, 0.5\}$ seconds.}
\label{fig:Academic Example}
\end{figure*}

We consider a system with dynamics

\begin{equation}\label{Academic Dynamics}
    \dot{x}(t) = f(x(t)) + G(x(t))u(t) = \begin{bmatrix} 10 u_1 + (3 - x_2)u_2 \\ (2 - x_1)u_1 + 7u_2 \\ (2.5 + x_3)u_3 \end{bmatrix},
\end{equation} 

\noindent where

\begin{equation*}
    f(x(t)) = \begin{bmatrix} 0 \\ 0 \\ 0 \end{bmatrix}, \indent G(x(t)) = \begin{bmatrix} 10 & 3 - x_2 & 0 \\ 2 - x_1 & 7 & 0 \\ 0 & 0 & 2.5 + x_3 \end{bmatrix},
\end{equation*}

\noindent with the primary interest of finding the reachable set of $x_1$ and $x_2$. Taking the Jacobian of $f(x)$ and $G(x)$ yields $L_f = L_G = 1$ to be acceptable Lipschitz constants. For simplicity, we set $x_0 = 0$. We remind the reader that $f$ and $G$ are assumed to be unknown, and only $L_f$ and $L_G$ are known.

Results from Figure \ref{fig:Academic Example}, showing the projection of the reachable sets to the first two coordinates, illustrate that (i) $\bar{\mathcal{R}}(T,0)$ and $\hat{\mathcal{R}}(T,0)$ are indeed underapproximations of the true reachable set, (ii) $\hat{\mathcal{R}}(T,0)$ produces a better underapproximation than $\bar{\mathcal{R}}(T,0)$, and (iii) the accuracy of the underapproximations increases as $T \to 0$. Phenomenon (i) validates the results from Theorems 3 and 4. The singular value decomposition of $G(0)$ yields $\sigma_1 \approx 11.43$, $\sigma_2 \approx 5.6$, and $\sigma_3 = 2.5$. Because of the large difference between $\sigma_1$ and $\sigma_3$, the true reachable set is not accurately represented by the ball $\bar{\mathcal{R}}(T,0)$. On the other hand, in Theorem 4 we derive control system \eqref{ODE2} using $\hat{\mathcal{V}}^\mathcal{G}_x$, which more accurately represents the complex shape of the GVS. Thus, the resulting set $\hat{\mathcal{R}}(T,0)$ is larger along singular vectors pertaining to larger singular values, contributing to result (ii).

Lastly, Figure \ref{fig:Academic Example} illustrates how the underapproximations $\bar{\mathcal{R}}(T,0)$ and $\hat{\mathcal{R}}(T,0)$ become asymptotically perfect as $T \to 0$. As $T \to \infty$, although comparatively worse, these underapproximations yield progressively larger reachable sets. In the next example, we apply this novel theory to decoupled quadrocopter dynamics to help safely land a damaged UAV. In this case, we will show that, in contrast to the above example, the polygonal approximation obtained by Theorem 4 produces worse results than the ball approximation derived from Theorem 3.

\subsection{Decoupled Quadrocopter Dymamics}

\begin{figure*}[h]
\centering
\subfloat{%
       \includegraphics[width=0.32\linewidth]{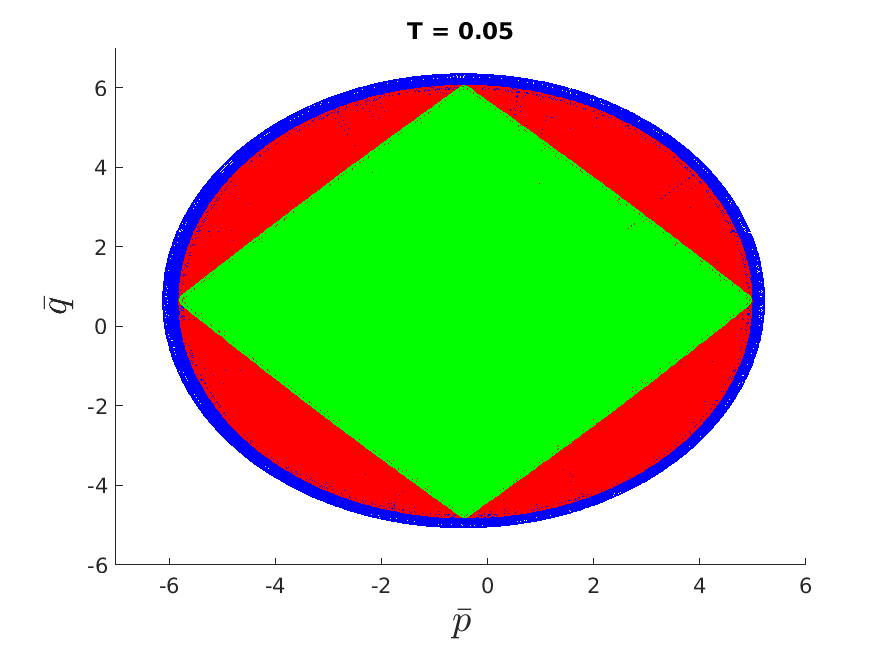}}
    \hfill
\captionsetup{justification=centering,margin=2cm}
\subfloat{%
       \includegraphics[width=0.32\linewidth]{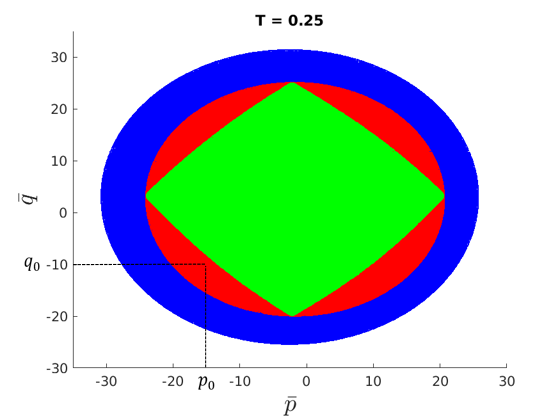}}
    \hfill
\captionsetup{justification=centering,margin=2cm}
\subfloat{%
       \includegraphics[width=0.32\linewidth]{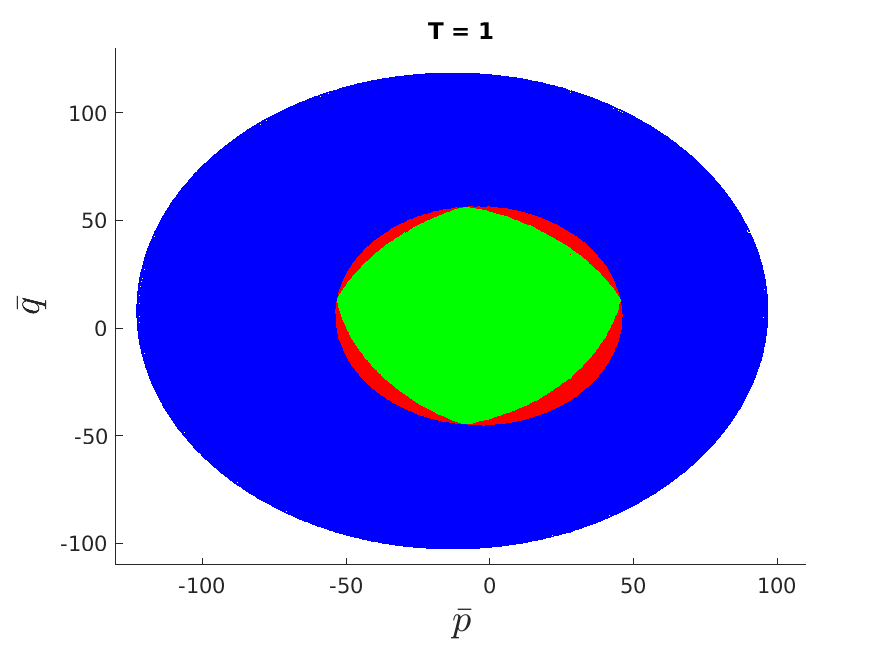}}
    \hfill
\captionsetup{justification=centering,margin=2cm}
\caption[font=small]{True reachable set (blue) with the underapproximations $\bar{\mathcal{R}}(T,0)$ (red) and $\hat{\mathcal{R}}(T,0)$ (green) numerically calculated for $T \in \{0.05, 0.25, 1\}$ seconds.}
\label{fig:Quadrocopter Example}
\end{figure*}

We show that the novel theory can be applied to a real system by example of quadrocopter dynamics. We consider the objective of adjusting the quadrocopter's pitch and roll velocities. For a safe landing, a UAV ideally needs these velocities to equal $0$ \citep{ghamry2016real}. Given the physical dimensions of a standard UAV, the yaw rate is inconsequential due to the symmetrical shape of the quadrocopter. We consider the scenario where a UAV collides with an obstacle, which would result in unwanted high velocity rotations. This scenario translates to a problem of reachability: given initial conditions of roll and pitch rates, we aim to determine if it is possible to reach $p = q = 0$ without knowing the system's dynamics after the collision. 

The dynamics of a standard UAV are modeled in \citep{beard2008quadrotor}. The model is comprised of a solid sphere with mass $M = 1\text{kg}$ and radius $R = 0.1\text{m}$, which represents the central frame; it is connected to four point masses $m = 0.1\text{kg}$, each representing one of four propellers at an equidistant length of $l = 0.5\text{m}$ away from the central sphere. The dynamics are shown below:

\begin{equation}\label{quadrotorDynamicsTrue}
    \begin{bmatrix} \dot{p} \\ \dot{q} \\ \dot{r} \end{bmatrix} = \begin{bmatrix} \frac{J_y - J_z}{J_x}qr \\ \frac{J_z - J_x}{J_y}pr \\ \frac{J_x - J_y}{J_z}pq \end{bmatrix} + \begin{bmatrix} \frac{1}{J_x} \tau_\phi \\ \frac{1}{J_y}\tau_\theta \\ \frac{1}{J_z}\tau_\psi \end{bmatrix},
\end{equation}

\noindent where states $p$, $q$, and $r$ correspond to velocities pertaining to the roll rate, pitch rate, and yaw rate, respectively, and

\begin{equation*}
    J_x = \frac{2MR^2}{5} + 2l^2m, \indent J_y = J_x, \indent J_z = \frac{2MR^2}{5} + 4l^2m.
\end{equation*}

Inputs $\tau_\phi$, $\tau_\theta$, and $\tau_\psi$ pertain to the applied torque that directly affects the roll, pitch, and yaw velocities. As $J_x = J_y$, the yaw rate can directly be changed by increasing the control action to $\tau_\psi$ without affecting the roll or pitch rates. Including additional states, such as the roll, pitch, and yaw angles and translational position, would move the model beyond the requirements of Assumption 1. Applying the theory to a dynamic structure such as this is a subject for future work.

We note the system dynamics in \eqref{quadrotorDynamicsTrue} are not globally Lipschitz continuous. However, as previously discussed, given the yaw rate can be directly changed by increasing control action to $\tau_\psi$ without affecting other states, we can trivially reduce $\dot{r}$ to $0$, causing $r$ to be some constant; we arbitrarily set $r_0 = \pi/2$. Let the initial conditions after collision be $p_0 = 15$, $q_0 = 10$ radians per second. Since novel theoretical results are derived under the assumption $x_0 = 0$, we perform a simple coordinate transformation; let $\bar{p} = p - 15$ and $\bar{q} = q - 10$ such that $\bar{p}_0 = \bar{q}_0 = 0$. Obviously, $\dot{\bar{p}} = \dot{p}$ and $\dot{\bar{q}} = \dot{q}$, thus resulting in the new system

\begin{equation}\label{quadrotorDynamicsTrue_Reduced}
    \begin{bmatrix} \dot{\bar{p}} \\ \dot{\bar{q}} \end{bmatrix} = \begin{bmatrix} \frac{\pi(J_y - J_z)}{2J_x}(\bar{q} + 10) \\ \frac{\pi(J_z - J_x)}{2J_y}(\bar{p} + 15) \end{bmatrix} + \begin{bmatrix} \frac{1}{J_x} \tau_\phi \\ \frac{1}{J_y}\tau_\theta \end{bmatrix},
\end{equation}

\noindent with initial conditions

\begin{equation*}
    f(0) = \begin{bmatrix} \frac{10\pi(J_y - J_z)}{2J_x} \\ \frac{15\pi(J_z - J_x)}{2J_y} \end{bmatrix} \approx \begin{bmatrix} -8.73 \\ 13.09 \end{bmatrix}, \indent  G(0) = \begin{bmatrix}  \frac{1}{J_x} & 0 \\ 0 & \frac{1}{J_y} \end{bmatrix},
\end{equation*}

\noindent such that $J_x = J_y = 0.009$, and $\|G(0)^\dagger\|^{-1} \approx 111.11$. The new dynamics \eqref{quadrotorDynamicsTrue_Reduced} are Lipschitz continuous; we overapproximate the Lipschitz bounds to be $L_f = 1$ and $L_G = 1$ to account for overly conservative knowledge about the rate of change of the dynamics. 

As in the previous sections, we bound inputs $u \in \mathbb{B}^2(0;1)$ and denote $x = [\bar{p}, \bar{q}]^T \in \mathbb{R}^2$ as the system states. The blue shape in Figure~\ref{fig:Quadrocopter Example} illustrates the true reachable set of the states pertaining to the roll and pitch velocities, in $(\bar{p},\bar{q})$ coordinates, plotted at $T \in \{0.05, 0.25, 1\}$ seconds. Given the structure of control systems \eqref{ODE1} and \eqref{ODE2}, it is trivial to see that if point $(-p_0,-q_0)$ is contained within either $\bar{\mathcal{R}}(T,0)$ or $\hat{\mathcal{R}}(T,0)$, then we know there exists a control $u$ which guarantees velocities $p$ and $q$ can be reduced to $0$, regardless of the true system dynamics.

Next, we apply novel theory to solve for $\bar{\mathcal{R}}(T,0)$, resulting in the set shown in red in Figure \ref{fig:Quadrocopter Example}, which is proven in Sections III and IV to be a guaranteed underapproximation of the GRS. We begin by calculating $\bar{\mathcal{V}}^\mathcal{G}_x = \mathbb{B}^2(f(0); \|G(0)^\dagger\|^{-1} - (L_f + L_G) \|x\|) = \mathbb{B}^2(f(0); 111.11 - 2\|x\|)$. Control system \eqref{ODE1} thus equals

\begin{equation} \label{UAV_ControlSystem1}
    \dot{x} = \begin{bmatrix} -8.73 \\ 13.09 \end{bmatrix} + \begin{bmatrix} 111.11 - 2\|x\| & 0 \\ 0 & 111.11 - 2\|x\| \\ \end{bmatrix} u
\end{equation}

\noindent such that $u \in \mathbb{B}^2(0;1)$. Its reachable set is $\bar{\mathcal{R}}(T,0)$.

Lastly, we consider control system \eqref{ODE2} to determine $\hat{\mathcal{R}}(T,0)$. Since $G(0)$ is invertible, we know from Corollary 2 that $\bar{\mathcal{V}}^\mathcal{G}_x \subseteq \bar{\bar{\mathcal{V}}}^\mathcal{G}_x$. Thus, when determining the diagonal entries of $\Lambda(\|x\|)$, we know from Theorem 4 that $\lambda_i(\|x\|) = \frac{g(\|x\|)}{\alpha(\|x\|)\|G(0)^\dagger \eta_i\| + \beta(\|x\|)}$. Since $G(0)$ is diagonal, we know the matrix of left singular vectors $U$ equals identity. System \eqref{ODE2} therefore equals

\begin{equation*}
    \dot{x} = \begin{bmatrix} -8.73 \\ 13.09 \end{bmatrix} + \begin{bmatrix} \lambda_1(\|x\|) & 0 \\ 0 & \lambda_2(\|x\|) \end{bmatrix} u
\end{equation*}

\noindent such that $\|u\|_1 \leq 1$.

The same methods for determining the previous two reachable sets numerically can be applied to find $\hat{\mathcal{R}}(T,0)$. Figure~\ref{fig:Quadrocopter Example} displays $\hat{\mathcal{R}}(T,0)$ for the roll and pitch velocities $\bar{p}$ and $\bar{q}$. Notice that $\|G(0)^\dagger \eta_i\|$ is identical for all $i$. Thus, $\lambda_i(\|x\|) = 111.11-2\|x\|$ for all $i$, which is identical to the diagonal terms in control system \eqref{UAV_ControlSystem1}. Since system \eqref{ODE2} considers $\|u\|_1 \leq 1$ and system \eqref{ODE1} considers $u \in \mathbb{B}^2(0;1)$, for this particular control system, $\hat{\mathcal{R}}(T,0) \subseteq \bar{\mathcal{R}}(T,0)$. The final result is a reachable set $\hat{\mathcal{R}}(T,0)$ denoted in green in Figure~\ref{fig:Quadrocopter Example}. It naturally resembles a polyhedron as $T \to 0$, with edges that curve as time increases.

According to the numerical results in Figure~\ref{fig:Quadrocopter Example}, roll and pitch velocities with initial conditions $p_0 = 15$, $q_0 = 10$ radians per second can be reduced to $0$ in no more than $0.25$ seconds because $(-15,-10)$ is contained in $\bar{\mathcal{R}}(0.25,0)$.  In fact, for any points that lie within $\bar{\mathcal{R}}(T,0)$, we can guarantee there exists a control signal which reaches these states within time $T$. For time $T = 0.05$, we see that neither the true reachable set, nor the underapproximations reach $(-15,-10)$, with $\bar{\mathcal{R}}(T,0)$ becoming asymptotically perfect as $T \to 0$. Conversely, at time $T = 1$, the theory provides larger, generally worse underapproximations where the true reachable set, along with $\bar{\mathcal{R}}(T,0)$ and $\hat{\mathcal{R}}(T,0)$, clearly include $(-15,-10)$.

Unlike in the first numerical example, set $\bar{\mathcal{R}}(T,0)$ is larger than $\hat{\mathcal{R}}(T,0)$. Naturally, in general, computing $\bar{\mathcal{R}}(T,0) \cup \hat{\mathcal{R}}(T,0)$ provides the best approximation of $\mathcal{R}^\mathcal{G}(T,0)$ available from our theory.
\section{Conclusion and Future Work}

This paper provides a novel approach to underapproximating the reachable set of a system with unknown dynamics. By assuming the nonlinear control-affine system structure and exploiting solely the knowledge of system dynamics at a single point and --- possibly conservative --- Lipschitz bounds on the rate of change, we are able to determine two underapproximations $\bar{\mathcal{R}}(T,x_0)$ and $\hat{\mathcal{R}}(T,x_0)$ that are guaranteed to be contained within the guaranteed reachable set $\mathcal{R}^{\mathcal{G}}(T,x_0)$. Both underapproximations rely on an intermediate approximation of the GRS by an ODI $\dot{x} \in \mathcal{V}^\mathcal{G}_x$, where the right hand side is a set of guaranteed velocities for the unknown nonlinear control-affine system. The two underapproximations differ by the shape of the right hand side set, i.e., the approximation of $\mathcal{V}^\mathcal{G}_x$ is either determined by balls $\bar{\mathcal{V}}^\mathcal{G}_x$ or a different shape $\bar{\bar{\mathcal{V}}}^\mathcal{G}_x$ that more closely resembles the shape of an intersection of an infinite many ellipsoids. 

A natural area of future work is to focus on creating a larger set underapproximating the GVS. One possibility is to consider approximations in norms other than the spectral norm considered in this paper. For example, potentially by bounding the perturbations of the unknown system's dynamics with the Frobenius norm instead of the spectral norm, we could utilize the Mirsky Inequality \citep{stewart1998perturbation} to produce new underapproximations of the GVS. However, there is currently no guarantee that such a bound would produce a more accurate underapproximation of the GVS. Another possibility includes the utilization of semi-infinite programming \citep{hettich1993semi, lopez2007semi} to determine the maximal extent of $\mathcal{V}^\mathcal{G}_x$ in every direction. Utilizing first- and second-order optimality conditions, the problem could potentially be reduced locally to one with finitely many constraints; similarly, given the convex structure of the GRS, additional optimization techniques such as duality may help further simplify the problem. 

Another approach to obtaining a larger underapproximation of the GVS is to increase the knowledge of the system dynamics. In other words, making additional assumptions on the structure of the dynamics could also help determine a larger underapproximation of the GVS by reducing the size of the set of systems consistent with prior knowledge about system dynamics. One option is to utilize knowledge from multiple system runs instead of solely dynamics at a single point. Expanding on theoretical results derived from this paper by incorporated additional knowledge consistent with a large class of systems could potentially result in substantial progress for the development of sophisticated safety critical systems.


\bibliographystyle{IEEEtran}
\bibliography{root.bib}

\begin{thebibliography}{10}
\providecommand{\url}[1]{#1}
\csname url@samestyle\endcsname
\providecommand{\newblock}{\relax}
\providecommand{\bibinfo}[2]{#2}
\providecommand{\BIBentrySTDinterwordspacing}{\spaceskip=0pt\relax}
\providecommand{\BIBentryALTinterwordstretchfactor}{4}
\providecommand{\BIBentryALTinterwordspacing}{\spaceskip=\fontdimen2\font plus
\BIBentryALTinterwordstretchfactor\fontdimen3\font minus
  \fontdimen4\font\relax}
\providecommand{\BIBforeignlanguage}[2]{{%
\expandafter\ifx\csname l@#1\endcsname\relax
\typeout{** WARNING: IEEEtran.bst: No hyphenation pattern has been}%
\typeout{** loaded for the language `#1'. Using the pattern for}%
\typeout{** the default language instead.}%
\else
\language=\csname l@#1\endcsname
\fi
#2}}
\providecommand{\BIBdecl}{\relax}
\BIBdecl

\bibitem{9304326}
M.~{Ornik}, ``Guaranteed reachability for systems with unknown dynamics,'' in
  \emph{59th IEEE Conference on Decision and Control}, 2020, pp. 2756--2761.

\bibitem{nguyen2008flight}
N.~Nguyen, K.~Krishnakumar, J.~Kaneshige, and P.~Nespeca, ``Flight dynamics and
  hybrid adaptive control of damaged aircraft,'' \emph{Journal of Guidance,
  Control, and Dynamics}, vol.~31, no.~3, pp. 751--764, 2008.

\bibitem{chowdhary2013guidance}
G.~Chowdhary, E.~N. Johnson, R.~Chandramohan, M.~S. Kimbrell, and A.~Calise,
  ``Guidance and control of airplanes under actuator failures and severe
  structural damage,'' \emph{Journal of Guidance, Control, and Dynamics},
  vol.~36, no.~4, pp. 1093--1104, 2013.

\bibitem{jourdan2010enhancing}
D.~Jourdan, M.~Piedmonte, V.~Gavrilets, D.~Vos, and J.~McCormick, ``Enhancing
  uav survivability through damage tolerant control,'' in \emph{AIAA Guidance,
  Navigation, and Control Conference}, 2010, pp. 7548--7569.

\bibitem{brockett1976nonlinear}
R.~W. Brockett, ``Nonlinear systems and differential geometry,''
  \emph{Proceedings of the IEEE}, vol.~64, no.~1, pp. 61--72, 1976.

\bibitem{isidori2013nonlinear}
A.~Isidori, \emph{Nonlinear Control Systems}.\hskip 1em plus 0.5em minus
  0.4em\relax Springer Science \& Business Media, 2013.

\bibitem{ornik2019control}
M.~Ornik, S.~Carr, A.~Israel, and U.~Topcu, ``Control-oriented learning on the
  fly,'' \emph{IEEE Transactions on Automatic Control}, vol.~65, no.~11, pp.
  4800--4807, 2019.

\bibitem{mitchell2003overapproximating}
I.~M. Mitchell and C.~J. Tomlin, ``Overapproximating reachable sets by
  hamilton-jacobi projections,'' \emph{Journal of Scientific Computing},
  vol.~19, no.~1, pp. 323--346, 2003.

\bibitem{filippova2018estimates}
T.~F. Filippova, ``Estimates of reachable sets of a nonlinear dynamical system
  with impulsive vector control and uncertainty,'' in \emph{14th International
  Conference: Stability and Oscillations of Nonlinear Control Systems}.\hskip
  1em plus 0.5em minus 0.4em\relax IEEE, 2018, pp. 1--4.

\bibitem{rungger2018accurate}
M.~Rungger and M.~Zamani, ``Accurate reachability analysis of uncertain
  nonlinear systems,'' in \emph{21st International Conference on Hybrid
  Systems: Computation and Control}, 2018, pp. 61--70.

\bibitem{mitchell2005time}
I.~M. Mitchell, A.~M. Bayen, and C.~J. Tomlin, ``A time-dependent
  hamilton-jacobi formulation of reachable sets for continuous dynamic games,''
  \emph{IEEE Transactions on Automatic Control}, vol.~50, no.~7, pp. 947--957,
  2005.

\bibitem{zhang2014reachable}
B.~Zhang, J.~Lam, and S.~Xu, ``Reachable set estimation and controller design
  for distributed delay systems with bounded disturbances,'' \emph{Journal of
  the Franklin Institute}, vol. 351, no.~6, pp. 3068--3088, 2014.

\bibitem{dullerud2013course}
G.~E. Dullerud and F.~Paganini, \emph{A Course in Robust Control Theory: A
  Convex Approach}.\hskip 1em plus 0.5em minus 0.4em\relax Springer Science \&
  Business Media, 2013.

\bibitem{ioannou2012robust}
P.~A. Ioannou and J.~Sun, \emph{Robust Adaptive Control}.\hskip 1em plus 0.5em
  minus 0.4em\relax Courier Corporation, 2012.

\bibitem{brunton2016discovering}
S.~L. Brunton, J.~L. Proctor, and J.~N. Kutz, ``Discovering governing equations
  from data by sparse identification of nonlinear dynamical systems,''
  \emph{Proceedings of the National Academy of Sciences}, vol. 113, no.~15, pp.
  3932--3937, 2016.

\bibitem{chen2018data}
Y.~Chen, H.~Peng, J.~Grizzle, and N.~Ozay, ``Data-driven computation of minimal
  robust control invariant set,'' in \emph{IEEE Conference on Decision and
  Control}.\hskip 1em plus 0.5em minus 0.4em\relax IEEE, 2018, pp. 4052--4058.

\bibitem{aubin2012differential}
J.-P. Aubin and A.~Cellina, \emph{Differential Inclusions: Set-Valued Maps and
  Viability Theory}.\hskip 1em plus 0.5em minus 0.4em\relax Springer Science \&
  Business Media, 2012.

\bibitem{bressan2007introduction}
A.~Bressan and B.~Piccoli, \emph{Introduction to the Mathematical Theory of
  Control}, 2007.

\bibitem{smirnov2002introduction}
G.~V. Smirnov, \emph{Introduction to the Theory of Differential
  Inclusions}.\hskip 1em plus 0.5em minus 0.4em\relax American Mathematical
  Society, 2002.

\bibitem{kurzhanski2000ellipsoidal}
A.~B. Kurzhanski and P.~Varaiya, ``Ellipsoidal techniques for reachability
  analysis,'' in \emph{International Workshop on Hybrid Systems: Computation
  and Control}.\hskip 1em plus 0.5em minus 0.4em\relax Springer, 2000, pp.
  202--214.

\bibitem{margaliot2007reachable}
M.~Margaliot, ``On the reachable set of nonlinear control systems with a
  nilpotent lie algebra,'' in \emph{European Control Conference}.\hskip 1em
  plus 0.5em minus 0.4em\relax IEEE, 2007, pp. 4261--4267.

\bibitem{vinter1980characterization}
R.~Vinter, ``A characterization of the reachable set for nonlinear control
  systems,'' \emph{SIAM Journal on Control and Optimization}, vol.~18, no.~6,
  pp. 599--610, 1980.

\bibitem{stewart1998perturbation}
G.~W. Stewart, ``Perturbation theory for the singular value decomposition,''
  Tech. Rep., 1998.

\bibitem{strang2016introduction}
G.~Strang, \emph{Introduction to Linear Algebra}.\hskip 1em plus 0.5em minus
  0.4em\relax Wellesley-Cambridge Press, 2016.

\bibitem{boyd2004convex}
S.~Boyd, S.~P. Boyd, and L.~Vandenberghe, \emph{Convex Optimization}.\hskip 1em
  plus 0.5em minus 0.4em\relax Cambridge University Press, 2004.

\bibitem{ross2013elementary}
K.~A. Ross, \emph{Elementary Analysis}.\hskip 1em plus 0.5em minus 0.4em\relax
  Springer, 2013.

\bibitem{stewart1977perturbation}
G.~W. Stewart, ``On the perturbation of pseudo-inverses, projections and linear
  least squares problems,'' \emph{SIAM Review}, vol.~19, no.~4, pp. 634--662,
  1977.

\bibitem{sontag1988controllability}
E.~D. Sontag, ``Controllability is harder to decide than accessibility,''
  \emph{SIAM Journal on Control and Optimization}, vol.~26, no.~5, pp.
  1106--1118, 1988.

\bibitem{chen2017high}
M.~Chen, \emph{High Dimensional Reachability Analysis: Addressing the Curse of
  Dimensionality in Formal Verification}.\hskip 1em plus 0.5em minus
  0.4em\relax University of California, Berkeley, 2017.

\bibitem{han2006reachability}
Z.~Han and B.~H. Krogh, ``Reachability analysis of nonlinear systems using
  trajectory piecewise linearized models,'' in \emph{American Control
  Conference}.\hskip 1em plus 0.5em minus 0.4em\relax IEEE, 2006, pp.
  1505--1510.

\bibitem{althoff2008reachability}
M.~Althoff, O.~Stursberg, and M.~Buss, ``Reachability analysis of nonlinear
  systems with uncertain parameters using conservative linearization,'' in
  \emph{47th IEEE Conference on Decision and Control}.\hskip 1em plus 0.5em
  minus 0.4em\relax IEEE, 2008, pp. 4042--4048.

\bibitem{chachuat2015set}
B.~Chachuat, B.~Houska, R.~Paulen, N.~Peri\'{c}, J.~Rajyaguru, and M.~E.
  Villanueva, ``Set-theoretic approaches in analysis, estimation and control of
  nonlinear systems,'' \emph{9th IFAC Symposium on Advanced Control of Chemical
  Processes}, vol.~48, no.~8, pp. 981--995, 2015.

\bibitem{ramdani2011computing}
N.~Ramdani and N.~S. Nedialkov, ``Computing reachable sets for uncertain
  nonlinear hybrid systems using interval constraint-propagation techniques,''
  \emph{Nonlinear Analysis: Hybrid Systems}, vol.~5, no.~2, pp. 149--162, 2011.

\bibitem{dormand1980family}
J.~R. Dormand and P.~J. Prince, ``A family of embedded runge-kutta formulae,''
  \emph{Journal of Computational and Applied Mathematics}, vol.~6, no.~1, pp.
  19--26, 1980.

\bibitem{althoff2018implementation}
M.~Althoff, D.~Grebenyuk, and N.~Kochdumper, ``Implementation of taylor models
  in {CORA} 2018,'' in \emph{5th International Workshop on Applied Verification
  for Continuous and Hybrid Systems}, 2018.

\bibitem{ghamry2016real}
K.~A. Ghamry, Y.~Dong, M.~A. Kamel, and Y.~Zhang, ``Real-time autonomous
  take-off, tracking and landing of uav on a moving ugv platform,'' in
  \emph{24th Mediterranean Conference on Control and Automation}.\hskip 1em
  plus 0.5em minus 0.4em\relax IEEE, 2016, pp. 1236--1241.

\bibitem{beard2008quadrotor}
R.~W. Beard, ``Quadrotor dynamics and control,'' 2008.

\bibitem{hettich1993semi}
R.~Hettich and K.~O. Kortanek, ``Semi-infinite programming: Theory, methods,
  and applications,'' \emph{SIAM Review}, vol.~35, no.~3, pp. 380--429, 1993.

\bibitem{lopez2007semi}
M.~Lopez and G.~Still, ``Semi-infinite programming,'' \emph{European Journal of
  Operational Research}, vol. 180, no.~2, pp. 491--518, 2007.

\end{thebibliography}

\end{document}